\documentclass{article}
\usepackage[utf8]{inputenc}

\input{def.tex}

\title{Conditional probability in Rényi spaces}
\author{Gunnar Taraldsen }
\date{\today}

\begin{document}

\maketitle

\begin{abstract}

  In 1933 Kolmogorov constructed a general theory that
  defines the modern concept of conditional probability.
  In 1955 Rényi fomulated a new axiomatic theory for probability motivated
  by the need to include unbounded measures.
  This note introduces a general concept of conditional probability in Rényi spaces.  
 
\end{abstract}


Kolmogorov is known for his axioms for probability,
but he himself gives much credit for this to french mathematicians, and Fréchet in particular.
\citet{KOLMOGOROV}
introduced, however, the definition of the
conditional expectation 
\be{1}
\pr^t (A) = \pr (A \st T = t)
\ee 
for the case where $\pr$ is a probability measure, 
$A$ is a positive measurable function, 
and $T: \Omega \into \Omega_T$ is a measurable mapping.
His definition extends the elementary definition, 
and covers in particular cases with $\pr(T=t) = 0$. 
\citet[p.v]{KOLMOGOROV} notes that
{\em especially the theory of conditional probabilities and conditional 
expectations} is an important novel contribution in his book.
Now, 86 years after the publication by \citet{KOLMOGOROV},
it is safe to say that Kolmogorov's note was most relevant.
The main result in this note is an extension of the Kolmogorov 
definition of $\pr^t$ to include the case where $\pr$ is a \Renyi state.

A \Renyi state $\pr$ is defined by
\be{2}
\pr (A \st B) = \frac{\mu(A B)}{\mu (B)}
\ee
where $\mu$ is a $\sigma$-finite measure on $\Omega$, 
$A$ is an event, 
and $B$ is an event with $0 < \mu (B) < \infty$.
An event is by definition a measurable set in $\Omega$.
The family $\cB$ of such elementary conditions $B$ 
obey the axioms
\begin{enumerate}
    \item $\emptyset \not\in \cB$
    \item $B,C \in \cB \imply B \cup C \in \cB$
    \item There exist $B_1, B_2, \ldots \in \cB$ with $\cup_i B_i = \Omega$. 
\end{enumerate}
A family $\cB$ of events that obeys these three axioms is
according to \citet[p.38, Def 2.2.1]{RENYI} a {\em bunch}.
A \Renyi state obeys the consistency relation
\be{3}
\pr (A \st B) = \frac{\pr (A B \st C)}{\pr (B \st C)},\;\;
\pr (B \st C) > 0
\ee
for all events $A$ and elementary conditions $B \subset C$.
A conditional probability space 
\citep[p.38, Def 2.2.2]{RENYI}
is a measurable space equipped 
with a consistent family of probability measures 
$\{\pr (\cdot \st B) \st B \in \cB\}$ where $\cB$ is a bunch.
It follows that a \Renyi state on 
$\Omega$ defines a conditional probability space. 

The structure theorem of \citet[p.40, Thm 2.2.1]{RENYI}
gives that any conditional probability space can be
represented by a $\sigma$-finite measure $\mu$ as in equation~(\ref{eq2}).
The original bunch $\cB'$ of the conditional probability space
is contained in the maximal bunch $\cB$ 
defined just after equation~(\ref{eq2}).
Equation~(\ref{eq2}) gives then also a maximal extension of the initial family of conditional probabilities since $\cB' \subset \cB$.
It follows also that a \Renyi state $\pr$ can be identified
with an equivalence class $[\mu] = \pr = \{c \mu \st c > 0\}$ of
$\sigma$-finite measures.
It will be convenient in the following to use the same symbol $\pr$ for both a \Renyi state and also for a $\sigma$-finite measure $\pr$ 
that represents the \Renyi state.
This is similar to the use of the symbol $f$ for a function and also 
a corresponding equivalence class $f$ of measurable functions. 
The study of spaces of \Renyi states is still in its infancy,
but we believe it will be important and interesting.
This is demonstrated by the study of convergence of Radon  \Renyi states by 
\citet{BiocheDruilhet16convergence},
and also - hopefully - by the constructive definition of conditional \Renyi states presented next.


We introduce now a definition of $\pr^t$ by defining $\pr^t (A \st B)$
for all positive measurable functions $A$ and all elementary conditions $B$.
It is, by the Radon-Nikodym theorem, uniquely determined by
requiring 
\be{4}
\int \phi (t) \pr^t (A \st B) \, \pr_T (dt \st B) =
\pr (\phi (T) A \st B)
\ee
to hold for all positive 
measurable functions $\phi$ using the definition
$\pr_T (C \st B) = \pr (T \in C \st B)$.
This is a natural generalization of the definition used by 
\citet[p.47, eq (1)]{KOLMOGOROV}.

We will next define a conditional \Renyi state $\pr^t$ as an
equivalence class $\pr^t = [\mu^t] = \{c (t) \mu^t \st \pr_T (c \le 0) = 0\}$
such that 
\be{5}
\pr^t (A B) = \pr^t (A \st B) \pr^t (B)
\ee
holds for all positive measurable functions 
$A$ and all elementary conditions $B$.
Choose a $\sigma$-finite measure $\nu$ that dominates $\pr_T$.
The proof of existence of $\nu$ is left to the reader.
A representative $\pr^t$ is now defined uniquely
by requiring 
\be{6}
\int \phi (t) \pr^t (A) \, \nu (dt) =
\pr (\phi (T) A)
\ee
to hold for all positive $\pr_T$ 
measurable functions $\phi$.
The result is a generalization of 
the disintegration of a measure relative to a pseudo-image 
as discussed by \citet[INT VI.45, published in 1959]{bourbaki03integration}.
\citet{TaraldsenTuftoLindqvist17improper} provide an alternative
route by constructing $\pr^t$ from $\pr^t (\cdot \st \cdot)$.

It remains only to prove that equation~(\ref{eq5}) holds.
Observe first that
$\pr_T (C \st B) = \pr ((T \in C) B) / \pr (B)$ and
$\pr ((T \in C) B) = \int_C \pr^t (B) \nu (dt)$ give
$\pr_T (dt \st B) = (\pr^t (B) / \pr (B)) \, \nu (dt)$.
Using this gives\\
$\int_C \pr^t (A B) \nu (dt) / \pr (B)$
$=$
$\pr ((T \in C) A \st B)$
$=$
$\int_C \pr^t (A \st B) \pr_T (dt \st B)$
$=$
$\int_C \pr^t (A \st B) \pr^t (B) \nu (dt) / \pr (B)$,
so
$$\int_C \pr^t (A B) \, \nu (dt) = \int_C \pr^t (A \st B) \pr^t (B) \, \nu (dt)$$
and equation~(\ref{eq5}) is proved since this holds for all events $C$ in
$\Omega_T$.
This proof is as given by 
\citet[Thm 1]{TaraldsenTuftoLindqvist18improper} who also give
a more detailed presentation with examples and motivation given by statistical 
inference including cases with improper priors and posteriors.
The results are key generalizations of the theory presented by \citet{TaraldsenLindqvist10ImproperPriors}.
It cannot be concluded, in the generality in this note, 
that a representative $\pr^t$ and $\pr^t (\cdot \st B)$ are measures for almost all $t$.
This was observed already by \citet[p.50,p.55]{KOLMOGOROV}.

\bibliographystyle{chicago}
\bibliography{main,bib}

\end{document}